\documentclass[12pt]{amsart}

\newtheorem{theorem}{Theorem}
\newtheorem{lemma}{Lemma}
\newtheorem{proposition}{Proposition}

\theoremstyle{definition}
\newtheorem{definition}{Definition}
\theoremstyle{remark}
\newtheorem{note}{Note}
\newtheorem{remark}{Remark}

\newcommand{\ltwo}{L^2({\mathbb R})}
\newcommand{\bbz}{\mathbb{Z}}
\newcommand{\bbrn}{\mathbb{R}^n}

\newcommand{\Bb}{\mathbb}

\begin{document}


\title{Robertson Type Theorems for Frames}
\author{Eric Weber}
\address{Department of Mathematics, Texas A\&M University, College Station, TX 77843-3368} 
\email{weber@math.tamu.edu}
\begin{abstract}
We extend Robertson's theorem to apply to frames generated by the action of a discrete, countable abelian unitary group.  Within this setup we use Stone's theorem and the theory of spectral multiplicity to analyze wandering frame collections.  Motivated by wavelet theory, we explicitly apply our results to the action of the integers given by translations on $\ltwo$.  This yields a new functional analytic method of constructing a wavelet from a multiresolution analysis.
\end{abstract}
\maketitle



\section{Introduction}

Frames have proven quite useful in a number of branches of mathematics and engineering; their study has taken quite an upturn in recent years.  Indeed, frames offer stable decomposition and reconstruction algorithms for analyzing signals, and are more flexible than orthonormal bases.  Frames also have important connections to multiresolution analyses and related structures.  This paper is motivated by those connections.  In this paper we shall consider frames for a Hilbert space $H$ generated by a group of unitary operators acting on a finite collection of vectors in $H$.

Let $G$ be a discrete, countable abelian group, and let $\pi: G \to B(H)$ be a unitary representation of $G$ on $H$.  Denote $\pi(g)$ by $\pi_g$.  We make the following definition:
\begin{definition}
A collection of vectors $W = \{w_1,\dots,w_n\}$ will be called a \emph{wandering frame collection} for $\pi_g$ if the collection $X = \{\pi_g w_i: g \in G, \ i=1,\dots,n\}$ is a frame for its closed linear span.  The collection $W$ will be called a \emph{complete wandering frame collection} if $X$ is a frame for $H$.  The collection $W$ will be called a \emph{wandering collection} or \emph{complete wandering collection} if $X$ is an orthonormal basis for its closed span or for $H$, respectively.
\end{definition}
Our purpose here is to describe the relationship between a wandering frame collection and a complete wandering frame collection.  Our results here can be considered as generalizations of Robertson's theorem \cite{R65}, which describes the relationship between a wandering subspace and a complete wandering subspace.  A wandering subspace $X \subset H$ of a unitary operator $V$ is a closed subspace such that $V^k(X) \perp V^l(X)$ for all $k \neq l$; a complete wandering subspace is such that the direct sum of all of the $V^k(X)$'s is $H$.  The statement of this theorem is:
\begin{theorem}[Robertson]
Let $X$ and $Y$ be wandering subspaces for a unitary operator $V$ such that:
\begin{enumerate}
\item[a.] $\sum_{k=-\infty}^{\infty} V^k(X) \subseteq \sum_{k=-\infty}^{\infty} V^k(Y)$,
\item[b.] dim(Y) $< \infty$.
\end{enumerate}
Then there exists a wandering subspace $\tilde{X}$ such that:
\begin{enumerate}
\item $X \subseteq \tilde{X}$,
\item $\sum_{k=-\infty}^{\infty} V^k(\tilde{X}) = \sum_{k=-\infty}^{\infty} V^k(Y)$
\end{enumerate}
\end{theorem}

The proof extends an orthonormal basis for $X$ to an orthonormal basis for $\tilde{X}$; this procedure is exact in the sense that there is a well defined number of wandering vectors needed to fill out $X$ to $\tilde{X}$, precisely $dim(Y) - dim(X)$.  This theorem is valid if we replace orthonormal bases with Riesz bases \cite{LTT95}; again the procedure is exact.

In our case with frames, we wish to, given a complete wandering frame collection, take a wandering frame collection and extend it so that it becomes a complete wandering frame collection.  However, there is no such well defined notion of size, only, in some sense a lower bound, which we shall make precise.  Indeed, one can take a frame for a Hilbert space and add as many 0 vectors as desired, and still have a frame.  As a result, our extension of Robertson's theorem will have several forms.

Our motivation arises from wavelet theory, in particular that of frame multiresolution analyses; for references see \cite{BL98a, P00a}.

Additionally, we shall be slightly more general than Robertson's theorem; we consider abelian groups of unitary operators instead of simply powers of a single unitary operator.  This is to accomodate the case of translations on $\Bb{R}^n$, which is a group that is not singly generated.

\section{Definitions and Main Results}

Robertson's theorem is fundamentally an extension theorem, i.e.~an extension of a wandering collection to a complete wandering collection.  Our main purpose is to show that a wandering frame collection can be extended to a complete wandering frame collection, or equivalently that a wandering frame collection always has a complementary wandering frame collection.

\begin{definition}
If $X = \{x_1, \dots, x_n\}$ and $Y = \{y_1, \dots, y_m\}$ are wandering frame collections for $\pi_g$ on $H$, we will say that $X$ and $Y$ are \emph{complementary} if $X \cup Y$ is a complete wandering frame collection for $\pi_g$ on $H$.
\end{definition}

\begin{note}
We do not require that the subspaces generated by $\{\pi_g x_i\}$ and $\{\pi_g y_j\}$ to be orthogonal.  We shall see that they may be chosen that way, however.
\end{note}

Our first result follows directly from a technical lemma given in section~\ref{S:CWFC}.
\begin{theorem} \label{T:main1}
Suppose that the representation $\pi_g$ admits a finite complete wandering frame collection $W$.  Suppose that $Y$ is a wandering frame collection for a subspace.  Then there exists a finite wandering frame collection that is complementary to $Y$.  This complementary collection may be chosen so that the resulting subspaces are orthogonal.
\end{theorem}

This result is not optimal in the sense that our complementary collection is in general bigger than necessary.  Our second result improves on the first in that our complementary collection $Y$ is smaller.
\begin{theorem} \label{T:main2}
Suppose that the representation $\pi_g$ admits a finite complete wandering frame collection $W$.  Suppose that $Y$ is a wandering frame collection for a subspace.  Then there exists an integer $k$, independent of the cardinality of $W$ and $Y$, such that:
\begin{enumerate}
\item there exists a finite wandering frame collection $X$ that is complementary to $Y$ and the cardinality of $X$ is $k$, and;
\item if $Z$ is complementary to $Y$, its cardinality is at least $k$.
\end{enumerate}
Moreover, the frame generated by $X$ can be taken to be normalized tight, regardless of whether the frames generated by $W$ and $Y$ are normalized tight.
\end{theorem}
This result can be considered as an extension of theorem 4 in \cite{HLPS}, and will also be proven in section~\ref{S:CWFC}.  Theorem~\ref{T:main2} has the drawback that the complementary collection $X$ is not always easily computable.  Therefore, the remainder of that section will discuss a ``middle ground'' procedure which can be computed but will reduce the cardinality of the complementary collection.  It is this procedure which generates a new algorithm for constructing a wavelet from a multiresolution analysis.  We outline the result here; the proof is in section~\ref{S:AW}.

\begin{remark}
Theorem~\ref{T:main1} is true even if the complete wandering frame collection $W$ is infinite, as will be evident in its proof.  However, in theorem~\ref{T:main2}, there may not exist such an integer.  Our motivation, multiresolution analyses, have a finite complete wandering frame collection except in extreme cases, which we will not consider here.
\end{remark}

Define unitary operators $D$, $T$ on $\ltwo$ as: $Df(x) = \sqrt{2} f(2x)$ and $Tf(x) = f(x - 1)$.  A wavelet is a function $\psi \in \ltwo$ such that the collection $\{D^n T^l \psi: n,\ l \in \Bb{Z} \}$ is an orthonormal basis for $\ltwo$.  A \emph{Multiresolution Analysis}, or MRA is a sequence of closed subspaces in $\ltwo$ that satisfy the following four conditions:
\begin{enumerate}
\item $V_j \subset V_{j+1}$,
\item $DV_j = V_{j+1}$,
\item $\cap_{j \in \bbz} V_j = \{0\}$ and $\cup_{j \in \bbz} V_j$ has dense span in $\ltwo$,
\item there exists a \emph{scaling function} $\phi \in V_0$, i.e. the collection $\{ T^l \phi: l \in \Bb{Z} \}$ is an orthonormal basis for $V_0$.
\end{enumerate}
By number 1.~above, we can define a second sequence of subspaces $\{W_j\}$ given by $V_{j+1} = V_j \oplus W_j$.  A routine calculation shows that $\{D\phi, DT\phi\}$ is a wandering collection for the unitary operator $T$ on $V_1$.  Whence by Robertson's theorem there exists a function $\psi \in W_0$ such that $\{T^l \psi: l \in \Bb{Z} \}$ is an orthonormal basis of $W_0$.  By numbers 2.~and 3.~above, it follows that such a function $\psi$ is a wavelet.  The proof of Robertson's theorem however, does not shed light on how to explicitly calculate such a function.

Mallat \cite{Ma89} provides an explicit algorithm for computing a wavelet from a MRA.  Our results presented in section~\ref{S:CWFC} provide an additional technique for constructing an orthonormal wavelet from a MRA.  Here is the idea of the result we obtain.  Define the functions
\[ \phi_0 = P_{W_0} D \phi; \quad \phi_1 = P_{W_0} DT \phi. \]

\begin{theorem}
Suppose that $\phi$ is an orthonormal scaling function for a multiresolution analysis.  Then there exist functions $f,g \in L^2([0,1))$ such that, with $f,g$ extended periodically, the function $\hat{\psi} = f \hat{\phi}_0 + g \hat{\phi}_1$ is the Fourier transform of an orthonormal wavelet.
\end{theorem}

\section{Complementary Wandering Frame Collections} \label{S:CWFC}

A sequence $\{f_i\}_{i \in I}$ is a frame for a separable Hilbert space $H$ if there exist positive constants A, B such that, for any $x \in H$,
\[ A\|x\|^2 \leq \sum_{i \in I}|\langle x, f_i \rangle|^2 \leq B\|x\|^2. \]
The frame $\{f_i\}_{i \in I}$ is called \emph{tight} if $A=B$ and \emph{normalized tight} if $A=B=1$.  We begin with a technical lemma.
\begin{lemma}
Let $\{x_i\}$ be a frame for the Hilbert space $H$, let $K$ be a closed subspace, and let $P$ be the projection of $H$ onto $K$.  Then $\{Px_i\}$ is a frame for $K$.  In particular, if $\{x_i\}$ is a normalized tight frame for $H$, then $\{Px_i\}$ is a normalized tight frame for $K$.
\end{lemma}
The proof appears in \cite{HL00a}, but we include it here for completeness.
\begin{proof}
For $z \in K$, we have
\[ A\|z\|^2 \leq \sum_{i \in I} |\langle z, x_i \rangle|^2 = \sum_{i \in I} |\langle Pz, x_i \rangle|^2 = \sum_{i \in I} |\langle z, Px_i \rangle|^2 \leq B\|z\|^2, \]
from which the lemma follows.
\end{proof}

\begin{proof}[Proof of theorem~\ref{T:main1}.]
Let $W = \{w_1, \dots, w_n\}$ and $Y = \{y_1, \dots, y_r\}$.
Let $H$ be the Hilbert space generated by $\pi_g W$, and $K$ be the subspace generated by $\pi_g Y$.  Let $P$ be the projection onto $K^\perp$.

Let $K$ be the subspace generated by $Y$, and let $P$ be the projection onto $K^{\perp}$.  Let $x_i = Pw_i$.  It follows from the lemma and the fact that $P$ commutes with the representation that $X = \{x_1, \dots x_n \}$ satisfies the statement.
\end{proof}

We now turn to proving theorem~\ref{T:main2}.  The proof relies on Stone's theorem for unitary representations of abelian groups, and the decomposition of projection valued measures, see \cite{BMM}.  However, we first require this technical result.
\begin{theorem} \label{T:rep}
A representation admits a finite complete wandering frame collection if and only if the representation is unitarily equivalent to a sub-representation some finite multiple of the regular representation.
\end{theorem}

A variation of this theorem appears in \cite[theorem 3.11]{HL00a}; there the statement is valid for any countable group but only for normalized tight frames.  We shall show, then, that any representation that admits a complete wandering frame collection also admits a complete wandering frame collection that generates a normalized tight frame.

Given a frame $\{f_n\}$ of a Hilbert space $H$, define an operator $S:H \to H$ given by $Sf = \sum_{n} \langle f, f_n \rangle f_n$; this is called the frame operator.  It is a positive, self-adjoint invertible operator, and the collection $\{S^{-1/2} f_n \}$ is a normalized tight frame for $H$.  See \cite{Cas00a}.
\begin{proof}
Suppose that $\{w_1,\dots,w_n\}$ is a complete wandering frame collection for $\pi_g$.  We first show that the frame operator $S$ commutes with the representation.  Let $f \in H$, $h \in G$ and compute:
\begin{align*}
S\pi_hf &= \sum_{i=1}^{n} \sum_{g \in G} \langle \pi_hf, \pi_g w_i \rangle \pi_g w_i \\
	&= \sum_{i=1}^{n} \sum_{g \in G} \langle f, \pi_{h^{-1}g} w_i \rangle \pi_g w_i \\
	&= \sum_{i=1}^{n} \sum_{g \in G} \langle f, \pi_g w_i \rangle \pi_{hg} w_i \\
	&= \pi_h \sum_{i=1}^{n} \sum_{g \in G} \langle f, \pi_g w_i \rangle \pi_g w_i \\
	&= \pi_h S f.
\end{align*}
Since $G$ is abelian, it follows that the $C^*$-algebra generated by $\pi(G)$ and $S$ is commutative, whence by the Spectral Theorem $S^{-1/2}$ commutes with $\pi(G)$.  Hence, the collection $\{S^{-1/2}w_1,\dots,S^{-1/2}w_n\}$ is a complete wandering normalized tight frame collection for $\pi$.
\end{proof}

\begin{proof}[Proof of theorem~\ref{T:main2}.]  
Our proof is based on a result in \cite{C}; we provide here the outline.  Let $K$ be the closed subspace spanned by $\{\pi(G)Y\}$.  By combining theorems~\ref{T:main1} and \ref{T:rep}, the representation of $G$ on $K^{\perp}$ is equivalent to a multiple of the regular representation.  Then by Stone's theorem and the theory of spectral multiplicity it follows that there exists a unitary operator
\[ U: K \to \oplus_{j=1}^{k} L^2(E_j, \mu), \]
where $\widehat{G} \supset E_1 \supset E_2 \supset \cdots \supset E_k$ and $\mu$ is the restriction of Haar measure to $E_1$, such that $U$ intertwines the projection valued measure on $H$ and the canonical projection valued measure on $\widehat{G}$.  Since $\{G\}$ forms an orthonormal basis of $L^2(\widehat{G}, \lambda)$, where $\lambda$ is Haar measure, the functions $g(\xi)\chi_{E_j}(\xi)$ form a normalized tight frame for $L^2(E_j,\mu)$.  It follows that the functions
\[ _gx_j(\xi) =  (0,\dots,\underbrace{g(\xi) \chi_{E_j}(\xi)}_{jth \ position}, \dots, 0) \] 
form a normalized tight frame for $UK$.  Hence, if we let $X = \{U^{-1} x_j \}$, then $X$ satisfies condition 1.~of the theorem.

To establish condition 2., note that the decomposition is unique (up to unitary equivalence), and that the summands are maximal cyclic subspaces, whence it follows that any cyclic collection (see definition~\ref{D:cyc} below) must have at least $k$ elements in it.
\end{proof}

\begin{note}
The unitary operator $U$ can be thought of as the Fourier transform on $K$.  Indeed, if $x \in K$, we will denote $Ux$ by $\hat{x}$.  We shall write $K \simeq \oplus_j L^2(E_j, \mu)$ when there exists a unitary operator $U$ that intertwines the projection valued measure and the canonical projection valued measure, as above.
\end{note}

We have now established our two main theorems.  However, as we stated earlier, the first proof yields (possibly) more vectors than we would want in practice, and the second theorem is not reasonably constructive.  So we wish to now demonstrate a technique for reducing the size of the collection from theorem~\ref{T:main1} by combining several vectors into one.

We first need to go back to the decomposition theorem as presented in the proof of theorem~\ref{T:main2}.  In the theory of the decomposition of projection valued measures, there exists a multiplicity function; in the case of the decomposition given in the proof above, the multiplicity function $m:\widehat{G} \to \Bb{N}$ is given by $m(\xi) = \sum_{j=1}^{n} \chi_{E_j}(\xi)$.  For our purposes, the multiplicity function will provide us a way of ``counting'' what parts of the regular representation we have.  The following lemma follows from the theory.

\begin{lemma} \label{L:mult}
Suppose $K = \overline{span}\{\pi_g x \}$, whence the representation on $K$ is cyclic.  Then the associated multiplicity function is given by $\chi_E$ where $E = \{\xi \in \widehat{G}: |\hat{x}(\xi)| > 0\}$.  If $X = \{x_1,\dots,x_k\} \in K$, then the multiplicity function associated to the subspace generated by $X$ is $\chi_F$ where $F = \{\xi \in \widehat{G}: max_{i=1,\dots,k}(|\hat{x_i}(\xi)|) >0\}$.
\end{lemma}

\begin{proposition} \label{P:prop1}
Suppose that the representation $\pi_g$ is cyclic on $K$.  Then the collection $\{w_1,\dots,w_k\}$ is a complete wandering frame collection if and only if the following two conditions hold:
\begin{enumerate}
\item $K \simeq L^2(E, \lambda|_E)$,
\item there exists $A,B > 0$ such that for almost every $\xi \in E$,
\[ A \leq \sum_{i=1}^{k} |\hat{w}_i(\xi)|^2 \leq B. \]
\end{enumerate}
Moreover, the frame bounds are given by 
\[ A = \text{ess inf} \sum_{i=1}^{k} |\hat{w}_i(\xi)|^2, \quad B = \text{ess sup} \sum_{i=1}^{k} |\hat{w}_i(\xi)|^2. \]
\end{proposition}
\begin{proof}
We have established the equivalence of condition 1.~to the existence of a complete wandering frame collection (for a cyclic representation).  Hence, we shall only consider condition 2.  We first show the necessity of the upper bound in condition 2.~by contrapositive.  Suppose $B > 0$ is given and suppose that $\sum_{i=1}^{k} |\hat{w}_i(\xi)|^2 > B$ for some set $F \subset E$ of positive measure.  Then $\hat{w}_i\chi_F \in L^2(E, \lambda)$, and consider the following calculation:
\begin{align*}
\sum_{i = 1}^{k} \sum_{g \in G} |\langle \chi_F, \widehat{\pi_g} \hat{w}_i \rangle|^2 &= \sum_{i = 1}^{k} \left| \int_{E} \chi_F(\xi) \overline{g(\xi) \hat{w}_i(\xi)} d\lambda \right| \\
&= \sum_{i = 1}^{k} \|\chi_F \hat{w}_i\|^2,
\end{align*}
since $G$ forms an orthonormal basis of $L^2(\widehat{G}, \lambda)$.  We have:
\begin{align*}
\sum_{i = 1}^{k} \|\chi_F \hat{w}_i\|^2 &= \sum_{i = 1}^{k} \int_{\widehat{G}} |\chi_F(\xi) \hat{w}_i(\xi)|^2 d\lambda \\
&= \int_{\widehat{G}} |\chi_F(\xi)|^2 \sum_{i = 1}^{k}|\hat{w}_i(\xi)|^2 d\lambda \\
&> B \int_{\widehat{G}} |\chi_F(\xi)|^2 d\lambda = B \|\chi_F\|^2,
\end{align*}
whence $B$ cannot be an upper frame bound.  The necessity of the lower bound in 2.~can be shown by an analogous calculation.

Likewise, to establish the sufficiency, let $x \in K$ and consider:
\begin{align*}
\sum_{i = 1}^{k} \sum_{g \in G} |\langle x, \pi_g w_i \rangle|^2 &= \sum_{i = 1}^{k} \sum_{g \in G} |\langle \hat{x}, \widehat{\pi_g} \hat{w}_i \rangle |^2 \\
&= \sum_{i = 1}^{k} \int_{\widehat{G}} \hat{x}(\xi) \overline{g(\xi) \hat{w}_i(\xi)} d\lambda \\
&= \sum_{i = 1}^{k} \|\hat{x}\hat{w}_i \|^2
\end{align*}
since, by condition 2., $\hat{x}\hat{w}_i \in L^2(E)$.  Moreover, by a calculation similar to above, 
\[ A\|\hat{x}\|^2 \leq \sum_{i = 1}^{k} \|\hat{x}\hat{w}_i \|^2 \leq B\|\hat{x}\|^2. \]
Finally, the frame bounds follow from a calculation analogous to the first calculation above.
\end{proof}

The idea we present here is to ``fuse'' two vectors from the wandering frame collection into one.  This cannot always be done.  When doing so we have two requirements: the first is that the resulting collection is cyclic for the entire space, and the second is that the collection retains frame bounds.

\begin{definition} \label{D:cyc}
A collection $X = \{x_1,\dots,x_m\} \subset H$ is called a \emph{cyclic collection} if the collection $\{\pi_g x_i\}$ has dense span in $H$.  A cyclic collection $X$ will be called \emph{reducible} if, after an appropriate reordering of $X$, there exists a $y_1 \in H$ such that $\{x_1,\dots,x_{n-2},y_1\}$ is also a cyclic collection.  We shall say that the vectors $x_{n-1}$ and $x_n$ are \emph{fusable}.
\end{definition}

Each element of a cyclic collection generates a cyclic subspace, whence we have the following lemma.

\begin{lemma} \label{L:fuse}
A cyclic collection is reducible if and only if two of its vectors are elements of a common cyclic subspace.  Equivalently, two vectors of a cyclic collection are fusable if and only if they are elements of a common cyclic subspace.
\end{lemma}
\begin{proof}
The only if implication follows from the preceding remark.  Clearly if $\{w_1,\dots,w_n\}$ is such that $w_{n-1}, w_n$ are elements of the cyclic subspace generated by $x_1$, then the collection $\{w_1, \dots, w_{n-2}, x_1\}$ is also cyclic.
\end{proof}

\begin{theorem}
Suppose the complete wandering frame collection $W = \{w_1,\dots,w_n\}$ for $\pi_g$ is reducible (as a cyclic collection), with the vectors $w_{n-1}$ and $w_n$ fusable.  Then $x_1$ can be chosen such that the cyclic collection $\widetilde{W} = \{w_1,\dots,w_{n-2},x_1\}$ is a complete wandering frame collection.
\end{theorem}
\begin{proof}
By lemma~\ref{L:fuse}, $w_{n-1}$ and $w_n$ are in a common cyclic subspace $K$.  Since the projection of $W$ onto $K$ yields a wandering frame collection for $K$ and $K$ is cyclic, we have that $K \simeq L^2(E,\lambda)$ for some $E \subset \widehat{G}$.

We shall construct $x_1$ in the following manner: let $F_1 = \{\xi : |\hat{w}_{n-1}(\xi)| \geq |\hat{w}_{n}(\xi)| \}$ and let $F_2 = E \setminus F_1$.  Then define $\hat{x}_1 = \hat{w}_{n-1}\chi_{F_1} + \hat{w}_{n}\chi_{F_2}$.  To show that $\widetilde{W}$ generates a frame, we need show that there exist frame bounds $\tilde{A},\tilde{B}$.  Let $A,B$ be the frame bounds for $W$; we shall show that $\tilde{A} = \frac{A}{2}$ and $\tilde{B} = 2 B $ suffice.  First suppose that $y \in K$.  By proposition~\ref{P:prop1}, $|\hat{x}_1|^2 \leq B$, whence by the calculation above,
\begin{align*}
\sum_{g \in G} | \langle y, \pi_g x_1 \rangle |^2 &= \|\hat{y}\hat{x}_1\|^2 \\
&= \int_{\widehat{G}} |\hat{y}|^2 |\hat{x}_1|^2 d\lambda \\
\intertext{and} \\
\sum_{g \in G} | \langle y, \pi_g w_{n-1} \rangle |^2 + \sum_{g \in G} | \langle y, \pi_g w_{n} \rangle |^2 &= \|\hat{y}\hat{w}_{n-1}\|^2 + \|\hat{y}\hat{w}_{n}\|^2 \\
&=\int_{\widehat{G}} |\hat{y}|^2 (|\hat{w}_{n-1}|^2 + |\hat{w}_n|^2) d\lambda.
\end{align*}
By our construction of $x_1$,
\[ \frac{1}{2}(|\hat{w}_{n-1}|^2 +|\hat{w}_{n}|^2) \leq |\hat{x}_1|^2 \leq (|\hat{w}_{n-1}|^2 +|\hat{w}_{n}|^2), \]
whence,
\begin{multline*}
\frac{1}{2} \left(\sum_{g \in G} | \langle y, \pi_g w_{n-1} \rangle |^2 + \sum_{g \in G} | \langle y, \pi_g w_{n} \rangle |^2 \right) \\
\leq \sum_{g \in G} | \langle y, \pi_g x_1 \rangle |^2 \leq \sum_{g \in G} | \langle y, \pi_g w_{n-1} \rangle |^2 + \sum_{g \in G} | \langle y, \pi_g w_{n} \rangle |^2.\\
\end{multline*}
Now let $y \in H$, $P_K$ be the projection onto $K$ and compute
\begin{align*}
\sum_{i = 1}^{n-2} & \sum_{g \in G} |\langle y, \pi_g w_i \rangle|^2 + \sum_{g \in G} |\langle y, \pi_g x_1 \rangle|^2 = \sum_{i = 1}^{n-2} \sum_{g \in G} |\langle y, \pi_g w_i \rangle|^2 + \sum_{g \in G} |\langle P_Ky, \pi_g x_1 \rangle|^2 \\
&\geq \frac{1}{2} \left(\sum_{i = 1}^{n-2} \sum_{g \in G} |\langle y, \pi_g w_i \rangle|^2 + \sum_{g \in G} |\langle P_Ky, \pi_g w_{n-1} \rangle|^2 + \sum_{g \in G} |\langle P_Ky, \pi_g w_{n} \rangle|^2 \right)\\
&\geq \frac{1}{2}A \|y\|^2.
\end{align*}
A similar calculation shows that $\widetilde{W}$ also has an upper frame bound of $2B$.
\end{proof}
\section{Application of our Results to Wavelets} \label{S:AW}

We now apply our results to wavelet theory.  The discrete countable abelian group in question is the integers, and the representation is on $\ltwo$ given by $\pi_l = T^l$.  In particular, given a MRA $\{V_j\}$, we restrict the action of $T^l$ to the subspace $V_1$.  Recall that if $\phi$ is the scaling function in $V_0$, we define the functions $\phi_0 = P_{W_0}D\phi$ and $\phi_1 = P_{W_0}DT\phi$.

Robertson's theorem shows that the representation of the integers on the subspace $W_0$ is cyclic--indeed it is unitarily equivalent to the regular representation of the integers.  Additionally, theorem~\ref{T:main1} shows that the collection $\{\phi_0, \phi_1\}$ is a complete wandering normalized tight frame collection for $\{T^l\}$ on $W_0$.  However, only one wandering vector is required, so we wish to ``fuse'' $\phi_0$ and $\phi_1$ as above.

We shall normalize the Fourier transform on $\ltwo$; for $f \in L^1(\Bb{R})\cap L^2(\Bb{R})$
\[\hat{f}(\xi) = \int_{\Bb{R}} f(x) e^{-2\pi i x \xi} dx. \]
For $\phi \in \ltwo$, define
\[|\vec{\hat{\phi}}(\xi)|^2 = \sum_{l \in \Bb{Z}} |\hat{\phi}(\xi + l)|^2. \]

In the case of the representation of the integers on $\ltwo$ given by translations, the multiplicity theory gives that
\[ \ltwo \simeq L^2(\widehat{G}, l^2(\Bb{Z})), \]
where the unitary operator yielding the equivalence is precisely the Fourier transform.  Hence, the representation on $W_0$ is equivalent to $L^2(\widehat{G}, \Bb{C})$ (via a unitary operator $U$), and also equivalent to a subspace of $L^2(\widehat{G}, l^2(\Bb{Z}))$ (via the Fourier transform), where each ``fiber'' is one dimensional.  For $f \in W_0$ then, we have that
\[ |Uf(\xi)|^2 = \sum_{l \in \Bb{Z}} |\hat{f}(\xi + l)|^2. \]
Moreover, we shall associate $\widehat{G} = S^1$ to $[0,1)$ in the standard way.
\begin{lemma}
Suppose that $\phi_0$ and $\phi_1$ are as above.  Let $E = \{ \xi : |\vec{\hat{\phi}}_0 (\xi)|^2 \geq \frac{1}{2} \}$, and let $F = \{ \xi : |\vec{\hat{\phi}}_1 (\xi)|^2 > \frac{1}{2} \}$.  Then the measure of $E \cap F$ is $0$, and $E \cup F = \Bb{R}$ up to a set of measure 0.
\end{lemma}
\begin{proof}
First note that both $E$ and $F$ are 1-periodic, i.e.~if $x \in E$, then $x+1 \in E$, whence we shall consider $E$ and $F$ restricted to $[0,1)$.  By theorem~\ref{T:main1}, the frame generated by $\{\phi_0, \phi_1\}$ is a normalized tight frame.  Since the representation on $W_0$ is equivalent to the regular representation, it's multiplicity function is identically 1 on $[0,1)$, and there exists a unitary operator $U: W_0 \to L^2(\widehat{G})$.  By our remarks above, 
\[ |U\phi_0(\xi)|^2 + |U\phi_1(\xi)|^2 = |\vec{\hat{\phi}}_0 (\xi)|^2 + |\vec{\hat{\phi}}_1 (\xi)|^2. \]
Since $\phi_0$ and $\phi_1$ together generate $W_0$, by lemma~\ref{L:mult}, $|\vec{\hat{\phi}}_0 (\xi)|^2 + |\vec{\hat{\phi}}_1 (\xi)|^2 > 0$ for almost all $\xi \in [0,1)$.  Finally, by proposition~\ref{P:prop1}, it follows that $|\vec{\hat{\phi}}_0 (\xi)|^2 + |\vec{\hat{\phi}}_1 (\xi)|^2 = 1$ for almost all $\xi \in [0,1)$, from which the lemma follows.
\end{proof}

\begin{theorem}
Let $\phi_0$ and $\phi_1$ be as above.  Define functions $f$ and $g$ by 
\[ f(\xi) = \frac{\chi_E(\xi)}{|\vec{\hat{\phi}}_0 (\xi)|}, \quad
g(\xi) = \frac{\chi_F(\xi)}{|\vec{\hat{\phi}}_1 (\xi)|}. \]
The function defined by
\[ \hat{\psi}(\xi) = f(\xi)  \hat{\phi}_0 (\xi) + g(\xi) \hat{\phi}_1 (\xi) \]
is an orthonormal wavelet for $\ltwo$, associated with the original MRA, i.e. $\psi \in W_0$.
\end{theorem}
\begin{proof}
By the definition of a MRA, it suffices only to show that the translates of $\psi$ form an orthonormal basis of $W_0$.  First note that $\psi \in W_0$ since, by definition, the functions $f$ and $g$ are 1-periodic and $f,g \in L^2([0,1))$.  Furthermore, we have that $|\vec{\hat{\psi}}(\xi)|^2 = 1$ a.e.~$\xi$, whence the translates of $\psi$ are orthonormal.  It follows then that the closed span of the translates of $\psi$ is $W_0$.
\end{proof}

Define the operator $T_{1/2}$ on $\ltwo$ by $T_{1/2}f(x) = f(x - 1/2).$  We have a much stronger result if the subspace $V_0$ reduces $T_{1/2}$.  (See \cite{W00a} for a discussion of when this occurs.)
\begin{theorem}
Suppose that the subspace $V_0$ reduces $T_{1/2}$. Then the function given by $\psi = \sqrt{2}\phi_0$ is an orthonormal wavelet.
\end{theorem}
\begin{proof}
The subspace $V_0$ reduces $T_{1/2}$ if and only if $W_0$ reduces $T_{1/2}$.  If so, then we have
\[ T_{1/2} \phi_0 = T_{1/2}P_{W_0}D\phi = P_{W_0}DT\phi = \phi_1. \]
Since $\hat{T}_{1/2}$ is a unitary multiplication operator, it follows that
\[ |\vec{\hat{\phi}}_0 (\xi)|^2 = |\vec{\hat{\phi}}_1 (\xi)|^2 = \frac{1}{2}. \]
Therefore, $\{T^l \phi_0:l \in \Bb{Z} \}$ forms a tight frame for $W_0$, with frame bound $\frac{1}{2}$.  Since the representation on $W_0$ is equivalent to the regular representation, $\{T^l \phi_0 \}$ forms an orthogonal set.  Whence, by normalizing $\phi_0$ by the factor of $\sqrt{2}$, we get an orthonormal basis of $W_0$, as required.
\end{proof}


\ifx\undefined\bysame
\newcommand{\bysame}{\leavevmode\hbox to3em{\hrulefill}\,}
\fi

\end{document}